\documentclass[12pt]{amsart}

\usepackage{amsmath,amsfonts,amscd,amsthm,amssymb,latexsym}

\font\teneufm=eufm10 \font\seveneufm=eufm7 \font\fiveeufm=eufm5
\newfam\frakturfam
\textfont\frakturfam=\teneufm \scriptfont\frakturfam=\seveneufm
\scriptscriptfont\frakturfam=\fiveeufm

\newtheorem{pr}{Proposition}

\newtheorem{lm}{Lemma}
\newtheorem{theor}{Theorem}
\newtheorem{co}{Corollary}

\def\bee{\begin{eqnarray}}
	\def\bes{\begin{eqnarray*}}
		\def\eee{\end{eqnarray}}
	\def\ees{\end{eqnarray*}}

\def\a{\alpha}
\def\b{\beta}

\def\Proof{{\sl Proof.}\ }

\title{Some examples of nonassociative coalgebras and supercoalgebras }

\begin{document}

\date{}
\maketitle

\begin{center}

{{\bf Daniyar Kozybaev}\footnote{
		Department of Mathematics,
		Eurasian National University, Nur-Sultan, Kazakhstan,
		e-mail:{\em kozybayev@gmail.com}},
	 {\bf Ualbai Umirbaev}\footnote{Department of Mathematics,
		Wayne State University,
		Detroit, MI 48202, USA; Department of Mathematics, 
		Al-Farabi Kazakh National University, Almaty, 050040, Kazakhstan; 
		and Institute of Mathematics and Mathematical Modeling, Almaty, 050010, Kazakhstan,
		e-mail: {\em umirbaev@wayne.edu}},
and {\bf Viktor Zhelyabin}\footnote{
		Institute of Mathematics of the SB of RAS, Novosibirsk, 630090, Russia,
		e-mail: {\em vicnic@math.nsc.ru}}, 
}
\end{center}

\begin{abstract} Locally finiteness of some varieties of nonassociative coalgebras is studied and the Gelfand-Dorfman construction for Novikov coalgebras and the Kantor construction for Jordan super-coalgebras are given. 
We give examples of a non-locally finite differential coalgebra, Novikov coalgebra, Lie coalgebra, Jordan super-coalgebra, and right-alternative coalgebra. The dual algebra of each of these examples satisfies very strong additional identities. We also constructed examples of an infinite dimensional simple differential coalgebra, Novikov coalgebra, Lie coalgebra, and Jordan super-coalgebra over a field of characteristic zero.   
  \end{abstract}



\section{Introduction}

\hspace*{\parindent}

The notion of coalgebra is dual to the notion of algebra. The theory of (co)associative coalgebras  has been developed for a long time within  the framework of the theory of Hopf algebras \cite{Swe69}. The study of coalgebras, bialgebras, and Hopf algebras received a new impetus when the term "quantum group", along with revolutionary new examples, was launched by V. Drinfeld in 1986  \cite{Dri87}.  Lie bialgebras, which are simultaneously Lie algebras and Lie coalgebras, were introduced as one of the most important notions of quantum group theory \cite{Dr83}. The study of Lie coalgebras, investigated earlier by W. Michaelis \cite{Michael80}, was intensified. It is well known that the dual of an associative coalgebra is an associative algebra and the dual of a Lie coalgebra is a Lie algebra. In 1994 J. Anquella, T. Cortes, and F. Montaner \cite{AnqCortMon94} called a coalgebra $C$  an $\mathfrak{M}$-coalgebra if the dual algebra $C^*$ belongs to the variety of algebras $\mathfrak{M}$. This allows to define alternative, Jordan, Malcev, left-symmetric, Novikov coalgebras, and so on. 

The Fundamental Theorem on Coalgebras asserts that every finitely  generated associative coalgebra over a field is finite-dimensional. 
An analogue of this result is true for alternative and Jordan coalgebras \cite{AnqCortMon94}, for structurable
coalgebras  \cite{Zhel95}, for Jordan copairs \cite{Zhel07}, for right alternative  Malcev admissible coalgebras and binary (-1,1)-coalgebras   \cite{SantMurShest21}.

 Lie coalgebras are not locally finite and the first example of an infinite dimensional finitely generated Lie coalgebra was given in \cite{Michael80}. In 1995 A. Slinko \cite{Slinko95} found some necessary and sufficient conditions for Lie coalgebras to be locally finite. A connection between Jordan and Lie (super)coalgebras, which is  an analog of the well known  Kantor-Koecher-Tits construction  for usual (super)algebras, was found in \cite{Zhel96,Zhel03}.  M.E.  Goncharov and V.N. Zhelyabin \cite{GonZh12, GonZh13} showed that every Malcev coalgebra embeds into a  Lie coalgebra with triality. Unlike Jordan coalgebras, Jordan super-coalgebras are not locally finite \cite{Zhel03}.
 
In 2000 D. Kozybaev \cite{Kozybaev00} constructed an example of a non-locally finite right-symmetric coalgebra and an example of a
non-locally finite right-alternative coalgebra. I. Shestakov reported \cite{SMS2} that the example of a right-alternative coalgebra given in \cite{Kozybaev00} is incorrect. This report attracted the attention of the authors to these old examples.
 First of all we noticed that the left-symmetric analogue of the non-locally finite right-symmetric coalgebra from \cite{Kozybaev00} is a Novikov coalgebra. Moreover, we noticed that the commutator coalgebra of this coalgebra is exactly the non-locally finite Lie coalgebra given by W. Michaelis \cite{Michael80}. In order to understand the nature of these examples we started to study codifferential coalgebras. 

To any associative and commutative differential algebra $A$ one relates the following three algebras: 

(1) a Novikov algebra obtained from $A$ by the  Gelfand-Dorfman construction; 

(2) a Lie algebra obtained as the commutator algebra of the Novikov algebra mentioned in (1); and 

(3) a Jordan superalgebra obtained from $A$ by the Kantor construction. 

The notion of coderivation allows us to define the notion of a (co)differential coalgebra.  We constructed a very easy example $C$ of a non-locally finite associative and commutative differential coalgebra. Using this example we constructed three examples of non-locally finite coalgebras. We define an analogue of the Gelfand-Dorfman construction for coalgebras  and using this we give an example of a non-locally finite Novikov coalgebra obtained from $C$, and the dual algebra of this coalgebra satisfies the identity $(xy)z=0$. This example coincides with Kozybaev's example mentioned above. Moreover, the commutator coalgebra of this Novikov coalgebra is exactly the non-locally finite Lie coalgebra given by Michaelis \cite{Michael80}. This Lie coalgebra is metabelian. Using $C$ and an analogue of the Kantor construction for super coalgebras, we give an example of a non-locally finite Jordan super-coalgebra. The dual of this super-coalgebra satisfies the super identities $xy=yx,xz=zx$, and $(zz_1)(z_2z_3)=0$ for even variables $x,y$ and for odd variables $z,z_1,z_2,z_3$. 

We repeated the same route starting from the simple differential algebra $(F[x],\partial=\frac{d}{d x})$ over a field $F$ of characteristic zero. The graded dual of this algebra is a simple  infinite dimensional differential coalgebra. Recall that a coalgebra is called simple if it does not have any proper subcoalgebras. Applying the Gelfand-Dorfman  construction to this coalgebra, we get an example of a simple infinite dimensional Novikov coalgebra. This coalgebra is the graded dual of the Novikov-Witt algebra  $\mathcal{L}_1$ \cite{KozUm16}.  
The commutator Lie coalgebra of this coalgebra is also simple and coincides with the graded dual of the Witt algebra $W_1$. 
And finally, using the Kantor construction, we constructed an example of a simple infinite dimensional Jordan super-coalgebra.

We also noticed that the example of a
non-locally finite right-alternative coalgebra from \cite{Kozybaev00} can be fixed only by interchanging the indexes $3n-2$ and $3n-1$ on lines 3 and 4 of the formula (12). The fixed example is given in Section 6. A much more complicated  example of a non-locally finite right-alternative coalgebra was recently given in \cite{SMS2}. 

The paper is organized as follows. In Section 2 we give some necessary terminology of coalgebras, notations, and general statements. In Section 3 we give examples of non-locally finite differential, Novikov, and Lie coalgebras. Simple infinite dimensional differential, Novikov, and Lie coalgebras are given in Section 4. Section 5 is devoted to Jordan super-coalgebras and Section 6 is devoted to right-alternative coalgebras.

\section{Coalgebras and coderivations}

\hspace*{\parindent}

Let $F$ be an arbitrary field. For any vector space   $V$ over $F$ denote by
 $$V^*={\rm Hom}_F(V,F)$$   its dual vector space, i.e., the vector space of all linear forms on  $V$. 

Denote by
$$V^{\otimes n}=\underbrace{V\otimes\ldots \otimes V}_n$$
 the $n$-th tensor power of the vector space $V$ over $F$. 

The map  
$$\rho: (V^*)^{\otimes n}\to (V^{\otimes n})^*$$ 
defined  by 
$$\rho(f_1\otimes\ldots\otimes f_n)(\sum_{i_1\ldots i_n}e_{i_1}\otimes\ldots\otimes  e_{i_n})=\sum_{i_1\ldots i_n}f_1(e_{i_1})\ldots f_n(e_{i_n})$$
is injective.  For this reason  we can  assume that $$(V^*)^{\otimes n}\subseteq   (V^{\otimes n})^*.$$ 

If  $\phi:V\to U$ is a linear map of vector spaces then the {\em transpose} $\phi^*: U^*\to V^*$ of $\phi$ is defined by 
$$\phi^*(u^*)(v)=u^*(\phi(v)), \  \  \ v\in V, u^*\in U^*.$$

\medskip

A vector space $C$ over  $F$ with a linear map 
$$\Delta: C\to C\otimes_FC$$
is called a {\em coalgebra}. The map $\Delta$ is called its 
 {\em comultiplication}. We often call the pair $(C,\Delta)$ a coalgebra in order to emphasize the comultiplication in question. 
For any  $a\in C$, using the Sweedler notation (see \cite{Swe69}), we write 
$$\Delta(a)=\sum_a a_{(1)}\otimes a_{(2)}.$$

If  $C$ is a coalgebra, then  
$$(fg)(a)=\rho(f\otimes g)(\Delta(a))=\sum_af(a_{(1)})g(a_{(2)}), \ \ \ \ f,g\in C^*, a\in C,$$
defines a product on $C^*$ and this product  turns $C^*$ into an algebra. Denote this product by 
\bes
m_{\Delta} : C^*\otimes C^*\to C^*. 
\ees
The algebra $C^*$ or $(C^*,m_{\Delta})$ is called the {\em dual algebra} of the coalgebra $(C,\Delta)$.

A coalgebra  $(C,\Delta)$ is called {\em coassociative}  if  
$$(\Delta\otimes \mathrm{id}-\mathrm{id}
 \otimes\Delta)\Delta=0,$$ i.e., for any $a\in C$ we have 
 $$\sum_a(\Delta(a_{(1)})\otimes a_{(2)}- a_{(1)}
 \otimes\Delta(a_{(2)})=0.$$

It is well known that a coalgebra $(C,\Delta)$  is coassociative if and only if its dual algebra  $C^*$ is associative. Moreover, a coalgebra $(C,\Delta)$  is a Lie coalgebra if and only if its dual  $C^*$ is a Lie algebra \cite{Michael80}. Following these results, the definition of coalgebras from any variety of algebras was given in \cite{AnqCortMon94}:
 
Let  $\mathfrak{M}$ be an arbitrary variety of algebras. A coalgebra $(C,\Delta)$  is called  an {\em $\mathfrak{M}$-coalgebra}  if its dual algebra  $C^*$ belongs to 
	$\mathfrak{M}$.

Let $V$ be a vector space  and let  
$\tau:V\otimes V\mapsto V\otimes V$ be the {\em ordinary flip}, i.e., a linear map with 
$\tau(x\otimes y)=y\otimes x$ for all $x,y\in V$. 

 A coalgebra  $(C,\Delta)$ is called {\em cocommutative} if 
$$\Delta=\tau\Delta,$$
 i.e., 
$$\sum_aa_{(1)}\otimes a_{(2)}=\sum_aa_{(2)}\otimes
 a_{(1)}$$ for any $a\in C$.

 Let $(C,\Delta)$ be an arbitrary coalgebra. A subspace $B$ of $C$ is called {\em subcoalgebra} of the coalgebra 
$(C,\Delta)$ if  $\Delta(B)\subseteq B\otimes B$. 

A subcoalgebra $B$ of a  coalgebra $(C,\Delta)$ is called {\em proper} if $B\neq \{0\}, C$. A coalgebra $(C,\Delta)$ without proper subcoalgebras is called {\em simple} \cite{Radford}.

It is well known \cite{AnqCortMon94} that $C$ admits a $C^*$-bimodule structure. The left and right actions of $C^*$ on $C$ are defined by 
$$\alpha\cdot a=\sum_{(a)}a_{(1)}\alpha(a_{(2)}), \ \ a\cdot \alpha=\sum_{(a)}\alpha(a_{(1)})a_{(2)}, \ \ \alpha\in C^*, a\in C.$$

Moreover,   a vector subspace 
$B$ of a coalgebra  $C$ is a subcoalgebra if and only if $B$ is a $C^*$-subbimodule  of $C$. Therefore the 
intersection of   any set of subcoalgebras of  $C$ is again  
a subcoalgebra.

Let  $S$ be a subset of  a coalgebra $C$. The smallest subcoalgebra $Coalg(S)$ of $C$ that contains  $S$ is called the subcoalgebra {\em generated by $S$}.  In other words,  $Coalg(S)$ is the $C^*$-subbimodule  of $C$ generated by  $S$.  If  $S$ is a finite set then  $Coalg(S)$ is called {\em finitely generated}.

A coalgebra $(C,\Delta)$ is called  {\em
	locally finite} if every    finitely generated subcoalgebra of  $ C$ is finite dimensional.

A linear map  $d:C\mapsto C$ is called a {\em coderivation} of the coalgebra  $(C,\Delta)$ if  
$$\Delta d=(d\otimes id+id\otimes d)\Delta,$$
i.e., 
$$\Delta(d(a))=\sum_a d(a_{(1)})\otimes a_{(2)}+a_{(1)}\otimes
d(a_{(2)}), \ \ a\in C. $$

A triple $(C,\Delta, d)$     is called a  {\em
	(co)differential coalgebra} if $(C,\Delta)$ is a coalgebra and $d$ is its coderivation. A subspace $B$ of a differential coalgebra $C$ is called a {\em subcoalgebra} if  $B$ is a subcoalgebra of $(C,\Delta)$ and $d(B)\subseteq B$, i.e., $B$ is a codifferentially closed subcoalgebra.

\begin{lm} \label{l1}
 Let  $d$ be a coderivation of the coalgebra $(C,\Delta)$. Then its transpose $d^*$  is a derivation of the dual algebra $C^*$, i.e., 
$$d^*(fg)=d^*(f)g+fd^*(g), \ \ f,g\in C^*.$$
\end{lm}
\Proof  Let $f,g\in C^*$, $a\in C$. Then
$$(d^*(fg))(a)=(fg)(d(a))=\rho(f\otimes
g)\Delta(d(a))=\sum_af(d(a_{(1)}))g(a_{(2)})+f(a_{(1)})g(d(a_{(2)}))$$
$$=\sum_a(d^*(f))(a_{(1)})g(a_{(2)})+f(a_{(1)})(d^*(g))(a_{(2)})=(d^*(f)g+fd^*(g))(a),$$
which proves the statement of the lemma. $\Box$

\begin{co} \label{c1}
 If $(C,\Delta, d)$ is an associative and commutative differential coalgebra then $(C^*,m_{\Delta}, d^*)$ is an associative and commutative differential algebra. 
\end{co}

  Let  $A$ be an algebra over $F$ with the multiplication  
$$m: A\otimes A\to A,$$ 
i.e., $m(a\otimes b)=ab$  for all  $a,b\in A$. Let 
$$m^*:A^*\mapsto (A\otimes A)^*$$
be the transpose of $m$. Unfortunately, the image $m^*(A^*)$ of $m^*$ does not always belong to $A^*\otimes A^*\equiv \rho(A^*\otimes A^*)\subseteq (A\otimes A)^* $. The structure of the dual coalgebra $(A^\circ,\Delta^\circ)$ is a little  complicated (see \cite{AnqCortMon94, Michael90}). 

Let 
\bes
A=\bigoplus_{i\in \mathbb{Z}} A_i, \ \ A_iA_j\subseteq A_{i+j}, 
\ees
be a $\mathbb{Z}$-graded algebra such that there exists an integer $m$ with $A_i=0$ for all $i<m$ and $A_i$ is finite dimensional for all $i\geq m$. In this case 
\bes
 (A\otimes A)_k=\bigoplus_{k=i+j} A_i\otimes A_j
\ees
is finite dimentional for all $k$ and 
\bes
 (A\otimes A)_k^*=\bigoplus_{k=i+j} \rho(A_i^*\otimes A_j^*). 
\ees
The coproduct $\Delta=\rho^{-1}m^*$ turns the graded space 
\bes
\bigoplus_{i\in \mathbb{Z}} A_i^*
\ees
into a coalgebra. This coalgebra is called the {\em graded dual} of the graded algebra $A$. 

\bigskip 

The following lemma is useful for studying subcoalgebras. 
\begin{lm}\label{l2}
	Let  $V$ be a vector space and let $W$ be a subspace of  $V$. Let  $a=\sum_{i=1}^na_i\otimes b_i$   and let the vectors  $a_1,\ldots,a_n$ be linearly independent. Assume  that  $a\in W\otimes W$. Then  $b_1, b_2,\ldots,b_n\in W$. 
\end{lm}
\Proof  Let $\alpha_1,\ldots,\alpha_n$ be a system of linear forms dual to $a_1,\ldots,a_n$, i.e., $\alpha_i(a_j)=\delta_{ij}$ for all $i,j$, where $\delta$ is the Kronecker delta function.

 Set $\phi = \alpha_1\otimes \mathrm{id} : V\otimes V\mapsto V$. Then 
 $$\phi(\sum e_i\otimes f_i)=\sum_i\alpha_1(e_i)f_i.$$ 
Obviously, $\phi(W\otimes W)\subseteq W$ and $\phi(a)=\sum_i\alpha_1(a_i)b_i=b_1$. Since  $a\in W\otimes W$ it follows that $\phi(a)\in W$. Consequently, $b_1\in W$. Similarly, we get $b_2,\ldots,b_n\in W$. $\Box$

 \section{Examples of non locally finite differential, Novikov, and Lie coalgebras}
 
 \hspace*{\parindent}

 An algebra  $A$ is called a {\em Novikov} algebra if it satisfaies the following identities: 
\bee\label{f1}
x(y z) - (x y) z = y (x z)-(y x) z,
\eee
\bee\label{f2}
 (x y) z =(x z) y.
\eee

Recall that an algebra satisfying the identity (\ref{f1}) is called {\em left-symmetric}.  Left-symmetric algebras are Lie-admissible, i.e., if $A$ is a left-symmetric algebra then $A$ with respect to the commutator $[x,y]:=xy-yx$ is a Lie algebra. This algebra is called the {\em commutator} algebra of $A$ and is denoted by $A^{(-)}$. 

The  identity (\ref{f1}) can be written as 
\bee\label{f3}
(x,y,z)=(y,x,z),
\eee 
where $(x,y,z):=(x y) z-x (y z)$ is the associator of elements $x,y,z$.

\bigskip

 {\bf The Gelfand-Dorfman construction \cite{GelDor79}.} Let $A$ be an associative and commutative algebra with a derivation  $d$. Define a new  multiplication $(\circ)$ on $A$ by 
$$x\circ y=xd(y),$$ where  $x,y\in A$.  Then $(A,\circ)$ is a Novikov algebra.

Moreover, the vector space $A$ with respect to the bracket 
\bes
[x,y]=xd(y)-yd(x)
\ees
is a Lie algebra. Obviously, $(A, [\cdot,\cdot])$ is the commutator algebra of the Novikov algebra $(A,\circ)$. 

Thus, any associative and commutative differential algebra $A$ generates the Novikov algebra $(A,\circ)$ and the Lie algebra $(A, [\cdot,\cdot])$. 
Over a field of characteristic zero every Novikov algebra can be embedded into a Novikov algebra $(A,\circ)$ for a suitable associative commutative differential algebra $A$ \cite{BCZ}. The class of Lie algebras embeddable into Lie algebras of the type $(A, [\cdot,\cdot])$ is not described yet.

 \begin{pr}\label{p1} A pair $(C,\Delta)$ is a Novikov coalgebra if and only if the following (co)identities hold:
\bee\label{f4}
(\Delta\otimes \mathrm{id}-\mathrm{id}\otimes\Delta)\Delta=(\tau\otimes \mathrm{id})(\Delta\otimes \mathrm{id}-\mathrm{id}\otimes\Delta)\Delta,
\eee
\bee\label{f5}
(\Delta\otimes \mathrm{id})\Delta=(\mathrm{id}\otimes \tau)(\Delta\otimes \mathrm{id})\Delta.
\eee
 \end{pr}
\Proof Let $(C,\Delta)$ be a coalgebra and let $\alpha,\beta,\gamma\in C^*$ and $c\in C$. Then 
$$(\alpha,\beta,\gamma)(c)=(\alpha\otimes \beta\otimes \gamma)((\Delta\otimes \mathrm{id}-\mathrm{id}\otimes\Delta)\Delta(c))$$
and 
$$(\beta,\alpha,\gamma)(c)=(\alpha\otimes \beta\otimes \gamma)((\tau\otimes \mathrm{id})(\Delta\otimes \mathrm{id}-\mathrm{id}\otimes\Delta)\Delta(c)).$$ 
Consequently, 
\bes
[(\alpha,\beta,\gamma)-(\beta,\alpha,\gamma)](c)\\
=(\alpha\otimes \beta\otimes \gamma)[(\Delta\otimes\mathrm{id}-\mathrm{id}\otimes\Delta)\Delta-(\tau\otimes \mathrm{id})(\Delta\otimes \mathrm{id}-\mathrm{id}\otimes\Delta)\Delta](c)). 
\ees
This implies that the identity (\ref{f3}) in $C^*$ is equivalent to the identity  (\ref{f4}) in $C$. 

Similarly, the identity (\ref{f2}) in $C^*$ is equivalent to the identity  (\ref{f5}) in $C$. $\Box$

\bigskip

{\bf The Gelfand-Dorfman construction for coalgebras.}  Let 
$(C,\Delta, d)$  be an associative and commutative differential coalgebra. Define on the space  $C$ a new comultiplication  $\Delta_N$  by   $$\Delta_N=(\mathrm{id}\otimes d)\Delta.$$  This means  
$$\Delta_N(a)=\sum_aa_{(1)}\otimes d(a_{(2)}) $$ for any  $a\in C$. Set also 
\bes
\Delta_L=\Delta_N^{(-)}=(1-\tau)\Delta_N, 
\ees
i.e., 
$$\Delta_N(a)=\sum_a(a_{(1)}\otimes d(a_{(2)})- d(a_{(2)})\otimes a_{(1)})$$ 
for any  $a\in C$. 

\begin{pr}\label{p2} $(1)$ 	The coalgebra  $(C,\Delta_N)$ is a Novikov coalgebra and the product in its dual algebra is defined by 
$$\alpha\circ\beta=\alpha
d^*(\beta), \ \ \ \ \alpha,\beta\in C^*.$$
$(2)$ The coalgebra  $(C,\Delta_L)$ is a Lie coalgebra and the bracket in its dual algebra is defined by 
$$[\alpha,\beta]=\alpha
d^*(\beta)-d^*(\alpha)\beta, \ \ \ \alpha,\beta\in C^*.$$
\end{pr}
\Proof By Corollary \ref{c1}, the dual $(C^*,m_{\Delta}, d^*)$ of the diferential coalgebra $(C,\Delta, d)$ is a differential algebra. We have 
$$(fg)(a)=\rho(f\otimes g)(\Delta(a)),$$ where  $f,g\in
C^*, \ a\in C$. By the Gelfand-Dorfman construction, the algebra  $(C^*,\circ)$ is a Novikov algebra, where $f\circ g=fd ^*(g)$ for 
all $f,g\in C^*$. 	
	
On the  other hand,
	$$(f\circ
	g)(a)=(fd^*(g))(a)=\sum_af(a_{(1)})d^*(g)(a_{(2)})=\sum_af(a_{(1)})g(d(a_{(2)}))=\rho(f\otimes
	g)(\Delta_N(a)),$$
i.e., 
$$f\circ
	g=fd^*(g)=\rho(f\otimes
	g)\Delta_N.$$

 Hence  $(C^*, \circ)$ is the dual algebra of the coalgebra  $(C,\Delta_N)$. Since  $(C^*, \circ)$ is a Novikov algebra it follows that $(C,\Delta_N)$ is a Novikov coalgebra. 

The second statement of the lemma can be checked similarly. $\Box$

\bigskip

{\bf Example 1.} Let  $C$ be a vector space  with a  linear basis  
\bes
e,f_1, f_2,\ldots,f_n,\ldots.
\ees
Define a  comultiplication 
 $\Delta:C\to C\otimes C$ on $C$  by  
\bes
\Delta(e)=e\otimes e, \ \Delta(f_i)=f_i\otimes e+e\otimes f_i, \ i\geq 1.
\ees
Define also a linear map  $d:C\to C$ by 
\bes
d(e)=0,\, d(f_i)=f_{i+1},\, i\geq 1.
\ees

\begin{lm}\label{l3}
  The triple   $(C,\Delta, d)$   is an associative and commutative differential coalgebra.  
\end{lm}
\Proof Obviously $(C,\Delta)$ is cocommutative. Direct calculations give 
that 
\bes
(\Delta\otimes \mathrm{id}-\mathrm{id}\otimes \Delta)\Delta(e)=e\otimes e\otimes e-
e\otimes e\otimes e=0, 
\ees
\bes
(\Delta\otimes \mathrm{id}-\mathrm{id}\otimes \Delta)\Delta(f_i)=
f_i\otimes e\otimes e+e\otimes f_i\otimes e\\
+e\otimes e\otimes f_i-
f_i\otimes e\otimes e-e\otimes f_i\otimes e-e\otimes e\otimes f_i=0,
\ees
i.e., $(C,\Delta)$ is coassociative. 

We also have 
$$\Delta(d(e))=0=d(e)\otimes e+e\otimes d(e),$$
$$\Delta(d(f_i))=\Delta(f_{i+1})=f_{i+1}\otimes e+e\otimes f_{i+1}=$$
$$(d(f_{i})\otimes e+d(e)\otimes f_{i+1})+(f_{i+1}\otimes d(e)+e\otimes d(f_{i})),\, i\geq 1.$$
This means that $d$ is a coderivation of the coalgebra
$(C,\Delta)$. $\Box$

By Lemma \ref{l3}, $(C,\Delta, d)$   is an associative and commutative differential coalgebra. Consequently, $(C^*,m_{\Delta}, d^*)$ is an associative and commutative differential algebra. Following the tradition in the theory of ordinary differential algebras, we denote the derivative of $x\in C^*$ by $x'$, i.e., $x'=d^*(x)$. 

\begin{pr}\label{p3}
The differential coalgebra $(C,\Delta, d)$ is not locally finite and the differential algebra $(C^*,m_{\Delta}, d^*)$ satisfies the differential identity 
$$x'y'=0.$$ 
\end{pr}
\Proof Obviously the codifferential subcoalgebra $B$ of $C$ generated by $f_1$ contains $f_i$ for all $i$. Since 
\bes
\Delta(f_1)=f_1\otimes e+ e\otimes f_1\in B\otimes B
\ees
it follows that $e\in B$ by Lemma \ref{l2}. Therefore $B=C$ is infinite dimensional. 

Let $\a,\b\in C^*$ and $c\in C$. Then
\bes
(d^*(\a)d^*(\b))(c)= \sum_c \a(d(c_{(1)}))\b(d(c_{(2)})). 
\ees
Consequently 
\bes
(d^*(\a)d^*(\b))(e)=\a(d(e))\b(d(e))=0
\ees
and 
\bes
(d^*(\a)d^*(\b))(f_i)=\a(d(f_i))\b(d(e))+\a(d(e))\b(d(f_i))=0
\ees
since $d(e)=0$. 

 This means that $C^*$ satisfies the differential identity $x'y'=0$. 
$\Box$

\bigskip

{\bf Example 2.} Consider the comultiplication $\Delta_N=(\mathrm{id}\otimes d)\Delta$ on $C$. 
By Proposition \ref{p2},  $(C,\Delta_N) $ is a Novikov coalgebra. We have 
$$\Delta_N(e)=(\mathrm{id}\otimes d)\Delta(e) =e\otimes d(e) =0$$ and 
$$\Delta_N(f_i)=(\mathrm{id}\otimes d)\Delta(f_i) =e\otimes d(f_i) =e\otimes f_{i+1}$$
for all $i\geq 1$.

\begin{theor}\label{t1}
The Novikov coalgebra $(C,\Delta_N)$ is not locally finite and the Novikov algebra  $C^*$ satisfies the identity
 $$(xy)z=0.$$
\end{theor}
\Proof  Let $B$ be the subcoalgebra of $(C,\Delta_N)$ generated by $f_1$. If $f_i\in B$ then 
$$\Delta_N(f_i)=(\mathrm{id}\otimes d)\Delta(f_i)=e\otimes
f_{i+1}\in  B\otimes B$$ 
implies that $e,f_{i+1}\in B$ by Lemma \ref{l2}. Consequently, $B=C$ is infinite dimensional.

Let  $\alpha,\beta,\gamma\in C^*$. Then 
$$((\alpha\beta)\gamma)(e)=(\alpha\otimes \beta\otimes \gamma)(\Delta\otimes \mathrm{id} )   \Delta(e)=0,$$
 $$((\alpha\beta)\gamma)(f_i)=(\alpha\otimes \beta\otimes \gamma)(\Delta\otimes \mathrm{id})   \Delta(f_i)=(\alpha\otimes \beta\otimes \gamma)(\Delta(e)\otimes f_{i+1} ) =0,$$
since  $\Delta(e)=0$. Hence the algebra   $C^*$ satisfies the identity  $(xy)z=0.$ $\Box$

Example 2 is the left-symmetric analogue of the non-locally finite right-symmetric coalgebra from \cite{Kozybaev00}.

\bigskip

{\bf Example 3.}
Now consider the Lie coalgebra $(C,\Delta_L)$. Recall that $\Delta_L=\Delta_N^{(-)}$ and, consequently,  
$$\Delta_L(e)=(d\otimes id+id\otimes d)\Delta(e)=0,$$ $$\Delta_L(f_i)=(d\otimes id+id\otimes d)\Delta(f_i)=e\otimes f_{i+1}-f_{i+1}\otimes e$$
for all $i\geq 1$. 

\begin{theor}\label{t2}
The Lie coalgebra $(C,\Delta_L)$ is not locally finite and the Lie algebra  $C^*$ satisfies the identity
 $$[[x,y],[z,t]]=0.$$
\end{theor}
\Proof Let $B$ be the subcoalgebra of $(C,\Delta_L)$ generated by $f_1$. If $f_i\in B$ then 
 $$\Delta_L(f_i)=e\otimes f_{i+1}-f_{i+1}\otimes e\in B\otimes B$$
implies that $e,f_{i+1}\in B$ by Lemma \ref{l2}. Consequently, $B=C$ is infinite dimensional.

Let  $\alpha,\beta,\gamma, \delta \in C^*$. Then 
$$([\alpha,\beta])(e)=(\alpha\otimes \beta)  \Delta_L(e)=0$$
and
\bes
([[\alpha,\beta],[\gamma,\delta]])(f_i)=(\alpha\otimes \beta\otimes \gamma\otimes\delta)(\Delta_L\otimes \Delta_L)   \Delta_L(f_i)\\
=(\alpha\otimes \beta\otimes \gamma\otimes\delta)(\Delta_L(e)\otimes \Delta_L(f_{i+1})-\Delta_L(f_{i+1})\otimes \Delta_L(e)) =0
\ees
since  $\Delta_L(e)=0$. Therefore $[[\alpha,\beta],[\gamma,\delta]]=0$. 	  $\Box$

Example 3 is Michaelis's example of a non-locally finite Lie coalgebra from \cite{Michael80}.

 \section{Infinite dimensional coalgebras without proper subcoalgebras}
 
 \hspace*{\parindent}

Let $F$ be a field of characteristic zero and let $F[x]$ be the algebra of polynomials over $F$ in one variable $x$.
Then  $(F[x], m, \partial)$, where $m$ is the polynomial multiplication and $\partial=\frac{d}{d x}$, is a simple differential algebra. Consider the natural grading  
\bes
F[x]=F 1 \oplus F x \oplus \ldots \oplus F x^n \oplus\ldots. 
\ees

\bigskip

The following example is the graded dual of the differential algebra $(F[x], m, \partial)$. 

{\bf Example 4.} Let 
\bes
C=F x_0 \oplus F x_1 \oplus \ldots \oplus F x_n \oplus\ldots, 
\ees
where $x_i\in F[x]^*$ is defined by $x_i(x^j)=\delta_{ij}$ for all $i,j$. It is easy to check that the dual comultiplication  
\bes
\Delta=\rho^{-1} m^* : C\to C\otimes C
\ees
 is defined by 
\bes
\Delta(x_n)=\sum_{i=0}^n x_i\otimes x_{n-i} 
\ees
and the coderivation 
\bes
d=\partial^*: C\to C 
\ees
is defined by 
 \bes
d(x_n)=(n+1)x_{n+1}
\ees
for all $n\geq 0$. 

Consequently, $(C,\Delta, d)$ is  an associative and commutative differential coalgebra. 

\begin{pr}\label{p4} The differential coalgebra $(C, \Delta, d)$ is simple. 
\end{pr}
\Proof Let $B$ be a nonzero subcoalgebra of $C$ and let 
\bes
f=x_n+\sum_{i<n} \a_ix_i \in B. 
\ees
Then 
\bes
\Delta(f)=x_n\otimes 1+\sum_{i<n} x_i\otimes f_i\in  B, \  \ f_i \in C. 
\ees
By Lemma \ref{l2} this implies that $x_n \in B$. Since $B$ is differentially closed it follows that $x_i\in B$ for all $i\geq n$. Morever, 
\bes
\Delta(x_n)=\sum_{i=0}^n x_i\otimes x_{n-i}\in B
\ees
imples that $x_0,\ldots,x_{n-1}\in B$ by Lemma \ref{l2}. Consequently, $B=C$. 
$\Box$ 

Applying the Gelfand-Dorfman construction to the differential algebra $(F[x],m,\partial)$ we get a Novikov algebra $(F[x], \circ)$, where $f\circ g= f\partial(g)$. This Novikov algebra was denoted by $\mathcal{L}_1$ in \cite{KozUm16} and is the first algebra in the list of left-symmetric Witt algebras $\mathcal{L}_n$. Notice that the commutator algebra of $\mathcal{L}_n$ is the Witt algebra $W_n$  \cite{KozUm16}. Since $\mathcal{L}_1$ is a Novikov algebra, we call $\mathcal{L}_1$ the {\em Novikov-Witt} algebra. Recall that 
\bes
\mathcal{L}_1=F \partial \oplus F x\partial \oplus \ldots \oplus F x^n\partial \oplus\ldots
\ees
is the algebra of all derivations of $F[x]$ with respect to the product
\bes
x^i\partial\circ x^j\partial=jx^{i+j-1}\partial, \ \ i,j\geq 0. 
\ees

\bigskip

The following example is the graded dual of the Novikov-Witt algebra $\mathcal{L}_1$. 

{\bf Example 5.} Let $(C,\Delta, d)$ be the differential coalgebra from Example 4. Applying the Gelfand-Dorfman construction, we get a Novikov coalgebra $(C, \Delta_N)$, where 
\bes
\Delta_N=(\mathrm{id}\otimes d)\Delta, 
\ees
i.e., 
\bes
\Delta_N(x_n)=\sum_{i=0}^n (n-i+1) x_i\otimes x_{n-i+1}
\ees
 
\begin{theor}\label{t3}
The Novikov coalgebra $(C, \Delta_N)$ is simple. 
\end{theor}
\Proof Let $B$ be a nonzero subcoalgebra of $C$ and let 
\bes
f=x_n+\sum_{i<n} \a_ix_i \in B. 
\ees
Then 
\bes
\Delta_N(f)=(n+1) 1\otimes x_{n+1}+\sum_{i<n+1} f_i\otimes x_i\in B, \  \ f_i \in C. 
\ees
By Lemma \ref{l2} this implies that $x_{n+1} \in B$. 
Consequently, we can assume that $x_n\in B$ for some $n\geq 0$. Then 
\bes
\Delta_N(x_n)=\sum_{i=0}^n (n-i+1) x_i\otimes x_{n-i+1}\in B
\ees
implies that $x_0,\ldots,x_{n+1}\in B$ by Lemma \ref{l2}.
 Consequently, $B=C$. $\Box$

\bigskip

The following example is the graded dual of the Witt algebra $W_1$. 

{\bf Example 6.} Let $(C,\Delta, d)$ be the Novikov coalgebra from Example 5. Then the commutator coalgebra $(C, \Delta_L)$ is a Lie coalgebra, where 
\bes
\Delta_L=\Delta_N^{(-)}=(1-\tau)\Delta_N,
\ees
i.e., 
\bes
\Delta_L(x_n)=\sum_{i=0}^n (n-i+1) x_i\otimes x_{n-i+1}-(i+1)x_{i+1}\otimes x_{n-i}\\
=\sum_{i=0}^{n+1} (n+1-2i) x_i\otimes x_{n+1-i}. 
\ees
 
\begin{theor}\label{t4}
The Lie coalgebra $(C, \Delta_L)$ is simple. 
\end{theor}
\Proof Let $B$ be a nonzero subcoalgebra of $C$ and let 
\bes
f=x_n+\sum_{i<n} \a_ix_i \in B. 
\ees
Then, as in the proof of Theorem \ref{t3}, we get $x_{n+1}\in B$. Applying Lemma \ref{l2} to the inclusion $\Delta_L(x_n)\in B\otimes B$, we get $x_0\in B$. Consequently, $B=C$. 
 $\Box$
 
Notice that an example of an infinite dimensional Lie coalgebra without finite dimensional subcoagebras was constructed in   \cite{Michael90} and \cite{Nichols90}.

\section{Non locally finite Jordan supecoalgebras}

\hspace*{\parindent}

Let  $G$ be the Grassman algebra with identity. Then $G=G_0+G_1$ is a 
$Z_2$-graded algebra. Let $J=J_0+J_1$ be a $Z_2$-graded algebra. Then $G(J)=J_0\otimes G_0+J_1\otimes
G_1$ is a subalgebra of the algebra  $G\otimes J$. The subalgebra $G(J)$ is called {\it Grassman envelope
} of the algebra $J$.

An algebra $J$ is called  a {\it Jordan superalgebra}, if its Grassman envelope  $G(J)$ is a Jordan algebra, i.e., $G(J)$ satisfies the following identities: 
$$xy=yx,$$
$$(x^2
y)x=x^2(yx).$$

\bigskip

{\bf The Kantor construction \cite{Kantor}.} Let $A$ be an associative commutative algebra over $F$ with a derivation  $D$. Denote by 
$\overline{A}$ an isomorphic copy of the vector space
$A$ with an isomorphism  $a\mapsto \overline{a}$. On the direct sum of the vector spaces 
$$J(A,D) = A\oplus\overline{A}$$ 
define a product $(\cdot)$  by 
$$a \cdot  b = ab,\, a \cdot \overline{b} =
\overline{ab},\, \overline{a}\cdot b = \overline{ab},\,
\overline{a}\cdot \overline{ b} = aD(b)-D(a)b,$$ where $a, b\in A$ and $ab$ is the product of elements in   $A$. Then $J(A,D)$ is a Jordan superalgebra. The superalgebra $J(A,D)$ is called a superalgebra of  the {\it vector type}.

\bigskip

{\bf The Kantor construction for coalgebras.} Let
$(C,\Delta,d)$ be an associtive and commutative differential coalgebra. Let $\overline{C}$ be an isomorphic copy of the vector space $C$ with an isomorphism $c\mapsto
\overline{c}$. On the direct sum of vector spaces 
$$J(C,d)=C\oplus\overline{C}$$
 define a coproduct $\Delta_J$ by 
\bes
\Delta_J(c)=\sum_{(c)}c_{(1)}\otimes c_{(2)}+\overline{
c_{(1)}}\otimes \overline{d(c_{(2)})}-\overline{d(c_{(1)})}\otimes
\overline{c_{(2)}},\\
\Delta_J(\overline{c})=\sum_{(c)}\overline{c_{(1)}}\otimes
c_{(2)}+c_{(1)}\otimes \overline{c_{(2)}},
\ees
 where  $c\in C$ and 
$\Delta(c)=\sum_{(c)}c_{(1)}\otimes c_{(2)}$. 

\begin{pr}\label{p5} The coalgebra $(J(C,d),\Delta_J)$ is a Jordan supercoalgebra and its dual 
	 $J(C^*,d^*)$ is a Jordan superalgebra of the vector type.
\end{pr}
\Proof By Corollary \ref{c1}, $(C^*,m_{\Delta}, d^*)$ is an associative and commutative differential algebra. We have 
$J(C,d)^*=C^*+(\overline{C})^*$. The isomorphism of $C$ and $\overline{C}$ induces the isomorphism of $C^*$ and $(\overline{C})^*$. Under this isomorphism, for any $\alpha\in C^*$ there corresponds $\overline{\alpha}\in (\overline{C})^*$ such that 
$\alpha(c)=\overline{\alpha}(\overline{c})$ for any  $c\in C$.
Therefore,  we can write that  $J(C,d)^*=C^*+\overline{C^*}$.

Let  $\alpha,\beta\in C^*$, and $c\in C$. Denote by  $(\cdot)$ the multiplication of  the algebra $(J(C,d),\Delta_J)^*$. Then we have
$$(\alpha\cdot\beta)(c)=\rho(\alpha\otimes\beta) \Delta_J(c),\,
 (\alpha\cdot\overline{\beta})(\overline{c})=\rho(\alpha\otimes\overline{\beta}) \Delta_J(\overline
 {c}),\,$$ $$ (\overline{\alpha}\cdot\beta)(\overline{c})=\rho(\overline{\alpha}\otimes\beta)\Delta_J(\overline
 {c}),\,(\overline{\alpha}\cdot\overline{\beta})(c)=\rho(\overline{\alpha}\otimes\overline{\beta})\Delta_J(c). $$
From this we get
$$\alpha\cdot\beta=\alpha\beta,\,
 \alpha\cdot\overline{\beta}=\overline{\alpha\beta},\, \overline{\alpha}\beta=\overline{\alpha\beta},\,
 \overline{\alpha}\cdot\overline{\beta}=\alpha d^*(\beta)-d^*(\alpha)\beta,
 $$ where  $\alpha\beta$ is the product of elements in the dual algebra $(C^*,m_{\Delta}, d^*)$.

 Consequently, $(J(C,d),\Delta_J)^*=J(C^*,d^*)$. $\Box$

\bigskip

{\bf Example 7.} Let $(C,\Delta, d)$ be the differential coalgebra from Example 1 and let $(J(C,d),\Delta_J)$ be the Jordan supercoalgebra obtained from $(C,\Delta, d)$ by the Kantor construction for coalgebras. Notice that 
$$\Delta_J(e)=e\otimes e, \ \ \Delta_J(f_i)=e\otimes f_i+f_i\otimes
e+\overline{e}\otimes
\overline{f_{i+1}}-\overline{f_{i+1}}\otimes\overline{ e},$$
$$\Delta_J(\overline{e})=e\otimes \overline{e}+\overline{e}\otimes e, \ \ \Delta_J(\overline{f_i})=e\otimes
\overline{f_i}+\overline{f_i}\otimes e+\overline{e}\otimes
f_i+f_i\otimes\overline{e},$$
for all $i\geq 1$. 

\begin{theor}\label{t5} The Jordan super-coalgebra  $(J(C,d),\Delta_J)$
	is not locally finite and its dual superalgebra satisfies the super identities 
\bes
xy=yx, xz=zx, (z_1z_2)(z_3z_4)=0
\ees
for all even $x,y$ and odd $z,z_1,z_2,z_3,z_4$. 
\end{theor}
\Proof Let  $B$ be the subcoalgebra of $J(C,d)$ generated by  $\overline{f_1}$.  If $\overline{f_i}\in B$, then $\Delta_J(\overline{f_i})\in B\otimes B$ implies that $e,\overline{e},f_i\in B$ by Lemma \ref{l2}.  If $f_i\in B$,  then 
$\Delta_J(f_i)\in B\otimes B$ implies that $\overline{f_{i+1}}\in B$.  Consequently,  $B=J(C,d)$ is infinite dimensional. 

By Proposition \ref{p5}, the dual of the supercoalgebra $(J(C,d),\Delta_J)$ is the Jordan superalgebra $J(C^*,d^*)$ obtained from the differential algebra $(C^*,m_{\Delta}, d^*)$ by the Kantor construction. The first two identities of $J(C^*,d^*)$, mentioned in the lemma, directly follow from the Kantor construction since  $C^*$ is an associative and commutative algebra. Notice that the product of two odd elements $z_1,z_2$ from $J(C^*,d^*)$ belongs to the ideal   $d^*(C^*)C^*$. By Proposition \ref{p3}, we have $d^*(C^*)^2=0$. Consequently, $(z_1z_2)(z_3z_4)=0$ for all odd elements $z_1,z_2,z_3,z_4$. 
$\Box$

The following example is the graded dual of the simple Jordan superalgebra $J(F[x],\partial)$ obtained from the simple differential algebra $(F[x],\partial)$ by the Kantor construction. 

{\bf Example 8.} Let $(C,\Delta, d)$ be the codifferential coalgebra from Example 4 and let $(J(C,d),\Delta_J)$ be the Jordan super-coalgebra obtained from $(C,\Delta, d)$ by the Kantor construction for coalgebras. Direct calculations give that 
$$\Delta_J(x_n)= \sum_{i=0}^n x_i\otimes x_{n-i} + \sum_{i=0}^{n+1} 
(n+1-2i )\overline{x}_i\otimes \overline{x}_{n-i+1},$$
$$\Delta_J(\overline{x}_n)= \sum_{i=0}^n( \overline{x}_i\otimes x_{n-i}+ x_i\otimes \overline{x}_{n-i}),$$
for all $n\geq 0$. 

\begin{theor}\label{t6} 
The Jordan supercoalgebra $(J(C,d),\Delta_J)$ is simple. 	
\end{theor}
\Proof Let $B$ be a nonzero subcoalgebra of $(J(C,d),\Delta_J)$ and let 
\bes
f=x_n+\sum_{i<n} \a_ix_i +\overline{c}\in B, 
\ees
where $\overline{c}\in \overline{C}$. Then 
\bes
\Delta_J(f)=\sum_{i\geq 0} a_i\otimes x_i +\sum_{i\geq 0} b_i\otimes \overline{x}_i \in B\otimes B
\ees
and it is easy to check that $b_{n+1}\neq 0$. By Lemma \ref{l2}, $\overline{x}_{n+1}\in B$. Then $\Delta_J(\overline{x}_{n+1})\in B\otimes B$ implies that 
$x_0,\ldots,x_{n+1},\overline{x}_0,\ldots,\overline{x}_{n+1}\in B$. Consequently, $B=C$.

If 
\bes
f=\overline{x}_n+\sum_{i<n} \a_i\overline{x}_i \in B, 
\ees
then $\Delta_J(f)\in B\otimes B$ implies that $x_n,\overline{x}_n\in B$ by Lemma \ref{l2}. Continuing the same discussions, we get $B=C$. 
 $\Box$

Notice that an example of a Jordan super-coalgebra without finite dimensional subcoagebras was constructed in  \cite{Zhel05}

\section{Non locally finite right alternative coalgebra}

\hspace*{\parindent}

In this section we give a corrected version of the example of a non locally finite right-alternative coalgebra from \cite{Kozybaev00}. This example was constructed on the base of the example of a finitely generated metabelian right-alternative algebra that is not residually finite \cite{UU89AA}.

 An algebra $A$ is called {\em right-alternative} if it satisfies  the following identities:
\bee\label{g1}
(yx)x=yx^2,
\eee
\bee\label{g2}
((xy)z)y=x((yz)y). 
\eee

The identity (\ref{g1}) is called {\em right-alternativity} and can be written in terms of the associators as 
\bee\label{g3}
(y,x,x)=0.
\eee
Over fields of characteristic $\neq 2$ this identity also implies the {\em Moufang identity} (\ref{g2}) (see \cite{ZSSS}).

\begin{pr}\label{p6} 	Let $(A,\Delta )$ be a coalgebra over a field $F$ of characteristic $\neq 2$. The coalgebra  $(A,\Delta )$ is right-alternative if and only if the following identity holds:   
\bee\label{g4}
(\Delta\otimes
	\mathrm{id}-\mathrm{id}\otimes\Delta)\Delta+(\mathrm{id}\otimes \tau)(\Delta\otimes
	\mathrm{id}-\mathrm{id}\otimes\Delta)\Delta=0.
\eee
\end{pr}
\Proof 
Let $(C,\Delta)$ be a coalgebra and let $\alpha,\beta,\gamma\in C^*$ and $c\in C$. Then 
$$(\alpha,\beta,\gamma)(c)=(\alpha\otimes \beta\otimes \gamma)((\Delta\otimes \mathrm{id}-\mathrm{id}\otimes\Delta)\Delta(c))$$
and 
$$(\alpha,\gamma,\beta)(c)=(\alpha\otimes \beta\otimes \gamma)((\mathrm{id}\otimes \tau)(\Delta\otimes \mathrm{id}-\mathrm{id}\otimes\Delta)\Delta(c)).$$ 
Consequently, 
\bes
[(\alpha,\beta,\gamma)+(\alpha,\gamma,\beta)](c)\\
=(\alpha\otimes \beta\otimes \gamma)[(\Delta\otimes\mathrm{id}-\mathrm{id}\otimes\Delta)\Delta+(\mathrm{id}\otimes \tau)(\Delta\otimes \mathrm{id}-\mathrm{id}\otimes\Delta)\Delta](c). 
\ees
This implies that the identity  (\ref{g4}) in $C$ is equivalent to linearized version of the identity (\ref{g3}) in $C^*$. Consequently,  (\ref{g4}) is equivalent to (\ref{g3}) over fields of characteristic $\neq 2$.  
$\Box$

\bigskip

{\bf Example 9.} Let $A$ be a vector space with a linear basis 
$$e_1 ,e_2,f_1,f_2,\ldots,f_n,\ldots.$$
 Define a comultiplication $\Delta$ on the vector space  $A$ by 
$$\Delta(e_1)=0,\,\Delta(e_2)=0,$$
$$  \Delta(f_{3n-2})=e_1\otimes f_{3n},\,\Delta(f_{3n-1})=e_2\otimes f_{3n},$$ $$\Delta(f_{3n})=e_2\otimes f_{3n+1}-f_{3n+1}\otimes e_2-e_1\otimes f_{3n+2}+f_{3n+2}\otimes e_1 ,\,   n\geq 1.$$

\begin{theor}\label{t7} The coalgebra  $(A,\Delta)$ is a right-alternative non locally finite coalgebra. Moreover, the dual algebra $A^*$ satisfies the identities $(xy)(zt)=((xy)z)t=0$. 
\end{theor}
\Proof Let  $\Delta_{Ass}=(\Delta\otimes \mathrm{id}-\mathrm{id}\otimes
\Delta)\Delta$. If $f, g, h\in A^*$, then
 $$(f,g,h)(a)=(f\otimes g\otimes h)(\Delta_{ass}
(a))$$ 
for all $a\in A$. We have $(f, g,g)(e_i)=0$ since $\Delta(e_i)=0$ for $i=1,2$.

 Direct calculations give that 
\bes
\Delta_{Ass}(f_{3n-2})=(\Delta\otimes \mathrm{id}-\mathrm{id}\otimes \Delta)(e_1\otimes f_{3n})=-e_1\otimes \Delta(f_{3n})\\
=-e_1\otimes e_2\otimes f_{3n+1}+e_1\otimes  f_{3n+1}\otimes e_2+e_1\otimes e_1 \otimes f_{3n+2}-e_1\otimes f_{3n+2}\otimes e_1, 
\ees
\bes
\Delta_{Ass}(f_{3n-1})=-e_2\otimes \Delta(f_{3n})
=-e_2\otimes e_2\otimes f_{3n+1}\\
+e_2\otimes  f_{3n+1}\otimes e_2+e_2\otimes e_1\otimes f_{3n+2}-e_2\otimes f_{3n+2}\otimes e_1, 
\ees
and
\bes
\Delta_{Ass}(f_{3n})=(\Delta\otimes \mathrm{id}-\mathrm{id}\otimes \Delta)(e_2\otimes f_{3n+1}-f_{3n+1}\otimes e_2-e_1\otimes f_{3n+2}+f_{3n+2}\otimes e_1 )\\
=-\Delta(f_{3n+1})\otimes e_2+\Delta(f_{3n+2})\otimes e_1 -e_2\otimes \Delta(f_{3n+1})+e_1\otimes \Delta(f_{3n+2})\\
=-e_1\otimes f_{3(n+1)}\otimes e_2+e_2\otimes f_{3(n+1)}\otimes e_1 -e_2\otimes e_1\otimes f_{3(n+1)}+e_1\otimes e_2\otimes f_{3(n+1)}.
\ees

Consequently, 
\bes
(f,g,g)(f_{3n-2})=
-f(e_1)g(e_2)g(f_{3n+1})+
f(e_1)g(f_{3n+1})g(e_2)\\
+
f(e_1)g(e_1)g(f_{3n+2})-f(e_1)g(f_{3n+2})g(e_1)=0, 
\ees 
\bes
(f,g,g)(f_{3n-1})=
-f(e_2)g(e_2)g(f_{3n+1})+
f(e_2)g(f_{3n+1})g(e_2)\\
+
f(e_2)g(e_1)g(f_{3n+2})-f(e_2)g(f_{3n+2})g(e_1)=0, 
\ees
and  
\bes
(f,g,g)(f_{3n})=
-f(e_1)g(f_{3(n+1)})g(e_2)+f(e_2)g( f_{3(n+1)})g(e_1)\\
-f(e_2)g(e_1)g(f_{3(n+1)})+f(e_1)g(e_2)g(f_{3(n+1)})=0.
\ees

Consequently, $(f,g,g)(a)=0$  for all  $f,g\in A^*$ and $a\in A$,
i.e.,   $(f,g,g)=0$. This means that $A^*$ satisfies the identity (\ref{g3}). 

Set $\phi=(\Delta\otimes \mathrm{id}\otimes
\mathrm{id})(\Delta\otimes \mathrm{id})\Delta$. Direct calculations give that 
$$\phi(e_1)=\phi(e_2)=0,$$ 
$$\phi(f_{3n-2})=(\Delta\otimes \mathrm{id}\otimes \mathrm{id})(\Delta\otimes \mathrm{id})(e_1\otimes f_{3n})=(\Delta\otimes\mathrm{id}\otimes \mathrm{id})(\Delta(e_1)\otimes  f_{3n})=0,$$
$$ \phi(f_{3n-1})=(\Delta\otimes \mathrm{id}\otimes \mathrm{id})(\Delta\otimes \mathrm{id})(e_2\otimes f_{3n})=(\Delta\otimes \mathrm{id}\otimes \mathrm{id})(\Delta(e_2)\otimes  f_{3n})=0 ,$$
and 
\bes
\phi(f_{3n})=(\Delta\otimes \mathrm{id}\otimes \mathrm{id} )(-\Delta(f_{3n+1})\otimes e_2+\Delta(f_{3n+2})\otimes
e_1)=\\
(\Delta\otimes \mathrm{id}\otimes \mathrm{id} )(-e_1\otimes f_{3(n+1)}\otimes e_2+e_2\otimes f_{3(n+1)}\otimes e_1)=0
\ees
for all $n\geq 1$. 

Consequently, $(((fg)h)e)(a)=0$ for all  $f,g,h,e\in A^*$ and
$a\in A$, i. e., $((fg)h)e=0$.  This means that $A^*$ satisfies the identity  
 $((xy)z)t=0$.

Now set $\psi=(\mathrm{id} \otimes\Delta \otimes \mathrm{id})(\mathrm{id}\otimes
\Delta)\Delta$. Then 
$$\psi(e_1)=\psi(e_2)=0,$$ 
\bes
\psi(f_{3n-2})=(\mathrm{id} \otimes\Delta \otimes \mathrm{id})(\mathrm{id}\otimes \Delta)(e_1\otimes f_{3n})=(\mathrm{id} \otimes\Delta \otimes \mathrm{id})(e_1\otimes \Delta(f_{3n}))\\
=
-e_1\otimes\Delta(f_{3n+1})\otimes e_2+e_1\otimes \Delta(f_{3n+2})\otimes e_1 \\
=-e_1\otimes e_1\otimes f_{3(n+1)}\otimes e_2+e_1\otimes e_2\otimes f_{3(n+1)}\otimes e_1,
\ees
\bes
\psi(f_{3n-   1})=(\mathrm{id} \otimes\Delta \otimes \mathrm{id})(\mathrm{id}\otimes \Delta)(e_2\otimes f_{3n})=(\mathrm{id} \otimes\Delta \otimes \mathrm{id})(e_2\otimes \Delta(f_{3n}))\\
=-e_2\otimes\Delta(f_{3n+1})\otimes e_2+e_2\otimes \Delta(f_{3n+2})\otimes e_1\\
=-e_2\otimes e_1\otimes f_{3(n+1)}\otimes e_2+e_2\otimes e_2\otimes f_{3(n+1)}\otimes e_1, 
\ees
and 
\bes
\psi(f_{3n})=(\mathrm{id} \otimes\Delta \otimes \mathrm{id})(\mathrm{id}\otimes \Delta)(e_2\otimes f_{3n+1}- f_{3n+1}\otimes e_2-e_1\otimes f_{3n+2}+f_{3n+2}\otimes e_1)\\
=(\mathrm{id} \otimes\Delta \otimes \mathrm{id})(e_2\otimes \Delta(f_{3n+1})
-e_1\otimes\Delta( f_{3n+2}))\\
=(\mathrm{id} \otimes\Delta \otimes
\mathrm{id})(e_2\otimes e_1\otimes f_{3(n+1)}-e_1\otimes e_2\otimes
f_{3(n+1)})\\
=e_2\otimes\Delta (e_1)\otimes f_{3(n+1)}-e_1\otimes
\Delta(e_2)\otimes f_{3(n+1)}= 0
\ees
for all $n\geq 1$. 

Therefore, 
\bes
(f((gh)g))(f_{3n-2})=-f(e_1)g(e_1)h(f_{3(n+1)})g(
e_2)+f(e_1)g(e_2)h(f_{3(n+1)})g(e_1)=0,
\ees 
\bes
(f((gh)g))(f_{3n-1})=0
\ees
 for all $f,g,h\in A^*$.

This means $f((gh)g)=0$
 for all  $f,g,h\in A^*$. Together with the identity $((xy)z)t=0$, this proves that the Moufang identity (\ref{g2}) holds in $A^*$. 

It is easy to check that $(\Delta\otimes \Delta)\Delta(a)=0$ for all  
$a\in A$. Therefore  $A^*$ satisfies the identity
$(xy)(zt)=0$. 

Now we show that$(A,\Delta)$ is not locally finite. 
Let  $B$ a subalgebra of  $(A,\Delta)$ generated by $f_1,f_2$.  

Suppose that $f_i\in B$. If $i=3n-1$ or $i=3n-2$, then  $\Delta(f_i)\in B\otimes B$ implies that $e_1,e_2,f_{3n}\in B$
 by Lemma \ref{l2}. If $i=3n$, then we get $e_1,e_2,f_{3n+1},f_{3n+2}\in B$. This implies that 
 $B=C$ is infinite-dimensional. $\Box$

\bigskip

\begin{center}
	{\large Acknowledgments}
\end{center}

The first author is supported by the grant of the Ministry of Education and Science of the Republic of Kazakhstan (project  AP09261086) and the other two authors are supported by the Russian Science Foundation (project 21-11-00286).

\hspace*{\parindent}

\end{document}